\documentclass{amsart}
\usepackage{amssymb}  

\theoremstyle{plain} 
\newtheorem{theorem}{Theorem}[section]
\newtheorem{lemma}[theorem]{Lemma}

\newtheorem{proposition}[theorem]{Proposition}
\newtheorem{corollary}[theorem]{Corollary}

\newtheorem{problem}[theorem]{Problem}
\theoremstyle{remark} 

\newtheorem{remark}[theorem]{Remark}

\numberwithin{equation}{section}

\newcommand{\seclabel}[1]{\label{sec:#1}} 
\newcommand{\thmlabel}[1]{\label{thm:#1}} 
\newcommand{\lemlabel}[1]{\label{lem:#1}} 
\newcommand{\corlabel}[1]{\label{cor:#1}} 
\newcommand{\prplabel}[1]{\label{prp:#1}} 
\newcommand{\prblabel}[1]{\label{prb:#1}}   

\newcommand{\secref}[1]{\ref{sec:#1}} 
\newcommand{\thmref}[1]{\ref{thm:#1}} 
\newcommand{\lemref}[1]{\ref{lem:#1}} 
\newcommand{\corref}[1]{\ref{cor:#1}} 
\newcommand{\prpref}[1]{\ref{prp:#1}} 
\newcommand{\prbref}[1]{\ref{prb:#1}}   

\newcommand{\Inn}{\operatorname{Inn}}

\newcommand{\lcm}{\mathrm{lcm}}

\newcommand{\ARIF}{\textsc{ARIF}}
\newcommand{\RIF}{\textsc{RIF}}
\newcommand{\CRIF}{\textsc{CRIF}}
\newcommand{\Mfg}{\textsc{Mfg}}
\newcommand{\Clp}{\textsc{C}}
\newcommand{\Flex}{\textsc{Flex}}

\newcommand{\RAlt}{\textsc{RAlt}}
\newcommand{\LAlt}{\textsc{LAlt}}
\newcommand{\LIP}{\textsc{LIP}}
\newcommand{\RIP}{\textsc{RIP}}
\newcommand{\AAIP}{\textsc{AAIP}}

\newcommand{\inv}{^{-1}}                                        
\newcommand{\linv}{^{\lambda}}                          
\newcommand{\rinv}{^{\rho}}                                 
\newcommand{\setof}[2]{\{ #1 \,|\, #2 \}}   
\newcommand{\sbl}[1]{\langle#1\rangle}          
\newcommand{\normal}{\trianglelefteq}               
\newcommand{\pc}[2]{#1_{\hskip -1pt(#2)}}       

\newcommand{\jeq}[1]{\stackrel{#1}{=}}          
\newcommand{\textjeq}[1]{\jeq{\text{#1}}}       

\newcommand{\Quasi}{\mathcal{Q}}                        
\newcommand{\Loop}{\mathcal{L}}                         
\newcommand{\TS}{\textbf{TS}$_0$\ }                 
\newcommand{\CWIP}{\textbf{CWIP}\ }                 
\newcommand{\Idem}{\mathcal{I}}                         
\newcommand{\Uni}{\mathcal{U}}                          

\author{Michael~K.~Kinyon}
\email{mkinyon@math.du.edu}
\urladdr{http://math.du.edu/\symbol{126}mkinyon}

\author{Petr~Vojt\v{e}chovsk\'y}
\email{petr@math.du.edu}
\urladdr{http://math.du.edu/\symbol{126}petr}

\address{Department of Mathematics \\
University of Denver \\
Denver, CO 80208 USA}

\title[Primary decomposition in loops]{Primary decompositions in varieties of commutative diassociative loops}

\keywords{commutative diassociative loop, commutative Moufang loop, commutative
C loop, $p$-primary component, RIF loop, ARIF loop, Steiner loop, Steiner
triple system, Moufang element, C element}

\subjclass{Primary: 20N05}

\begin{document}

\begin{abstract}
The decomposition theorem for torsion abelian groups holds analogously for
torsion commutative diassociative loops. With this theorem in
mind, we investigate commutative diassociative loops satisfying the
additional condition
(trivially satisfied in the abelian group case) that all $n$th powers are
central, for a fixed $n$. For $n=2$, we get precisely commutative $C$ loops.
For $n=3$, a prominent variety is that of commutative Moufang loops.

Many analogies between commutative C and Moufang loops have been noted in the
literature, often obtained by interchanging the role of the primes $2$ and $3$.
We show that the correct encompassing variety for these two classes of loops is
the variety of commutative RIF loops. In particular, when $Q$ is a commutative
RIF loop: all squares in $Q$ are Moufang elements, all cubes are $C$ elements,
Moufang elements of $Q$ form a normal subloop $M_0(Q)$ such that $Q/M_0(Q)$ is
a C loop of exponent $2$ (a Steiner loop), C elements of $L$ form a normal
subloop $C_0(Q)$ such that $Q/C_0(Q)$ is a Moufang loop of exponent $3$. Since
squares (resp. cubes) are central in commutative C (resp. Moufang) loops, it
follows that $Q$ modulo its center is of exponent $6$. Returning to the
decomposition theorem, we find that every torsion, commutative RIF loop is a
direct product of a C $2$-loop, a Moufang $3$-loop, and an abelian group
with each element of order prime to $6$.

We also discuss the definition of Moufang elements, and the quasigroups
associated with commutative RIF loops.
\end{abstract}

\maketitle

\section{Introduction}
\seclabel{intro}

A \emph{quasigroup} $(Q,\cdot)$ is a set $Q$ with a binary operation
$\cdot$ such that for each $a,b\in Q$, the equations $ax=b$, $ya=b$ have
unique solutions $x,y\in Q$, respectively. A \emph{loop} is a
quasigroup with a neutral element $1$, i.e., $1x=x1=x$ for every $x$.
Basic references for quasigroups and loops are \cite{Bruck,Pflugfelder}.

A loop is
\emph{power-associative} if every element generates a subgroup
(associative subloop), and
\emph{diassociative} if every two elements generate a subgroup. Powers $x^n$
are thus defined unambiguously in power-associative loops, and the order $|x|$ of
$x$ can be introduced in the usual way.

For a power-associative loop $Q$ and a prime $p$, the $p$-\emph{primary
component} $\pc{Q}{p}$ is the set of all torsion elements $x\in Q$ such that
$|x|$ is a power of $p$. A power-associative loop $Q$ is a $p$-\emph{loop} if
$Q=\pc{Q}{p}$.

A classical theorem of group theory states that every finitely generated
torsion abelian group is a direct product of its $p$-primary components. For
power-associative loops, a $p$-primary component need not even be a
subloop. On the other hand, Bruck and Paige observed without proof in
\cite{BruckPaige} that the decomposition theorem holds in the variety of
commutative diassociative loops. (We give a proof in \S\secref{general}.)

In this paper, we investigate the situation when additional equational
restrictions are imposed on the $p$-primary components of commutative
diassociative loops.

The condition that all $n$th powers (for a fixed $n$) are central is trivially
satisfied for commutative groups but not so for commutative diassociative
loops, since the center of a loop consist of all elements that commute
\emph{and} associate with all other elements.

In fact, the situation is fully understood only for the variety of
commutative diassociative loops with squares in the center; this coincides
with the variety of commutative C loops.
A loop is called a \emph{C loop} if it satisfies the identity
\[
x(y\cdot yz) = (xy\cdot y)z\,.  \tag{\Clp}
\]
C loops satisfying $x\cdot yx = xy\cdot x$, which include the
commutative ones, are diassociative \cite{KKP}.

The variety of commutative diassociative loops with cubes in the
center includes commutative Moufang loops.
A loop is called a \emph{Moufang loop} if
it satisfies any, and hence all, of the equivalent identities
\[
\tag{\Mfg}
\begin{array}{ccc}
x(yz\cdot x) = xy\cdot zx\,, &\qquad& (x\cdot yz)x = xy\cdot zx\,,\\
x(y\cdot xz) = (xy\cdot x)z\,, &\qquad& (zx\cdot y)x = z(x\cdot yx)\,.
\end{array}
\]
The diassociativity of Moufang loops is usually known as
Moufang's Theorem \cite{Bruck,Pflugfelder}.

Already for $n=3$ do we find that the variety of
commutative diassociative loops with central $n$th powers (for a fixed
$n$) is rather unwieldy, because
it properly contains the variety of commutative Moufang loops.
For instance, from the general construction of Hart and Kunen \cite{HK},
there exist nonMoufang, commutative diassociative loops of exponent $3$ and
order $27$.

Thus, although the decomposition theorem for the variety of commutative
diassociative loops with central $n$th powers is easy to prove
(see \S\secref{general}), it is not particularly useful, because this
variety is too broad. Ideally, we would like to be able to characterize
subvarieties of commutative diassociative loops whose $p$-primary components
satisfy certain prescribed (equational) conditions. In general, however, this
seems to be a difficult task.

In our previous work \cite{PhillipsVojtechovsky}, we observed many analogies
between commutative C loops and commutative Moufang loops, with $p=2$ playing a
prominent role in the C case and $p=3$ in the Moufang case. For instance,
as we have already noted,
squares of elements in a commutative C loop are central, while cubes of
elements in a commutative Moufang loop are central. In addition, a commutative
C loop is a direct product of an abelian group and a commutative C $2$-loop,
while a commutative Moufang loop is a direct product of an abelian group and a
commutative Moufang $3$-loop. The present work was in part motivated by our
desire to better understand this analogy.

It turns out that the behavior of commutative C and commutative Moufang loops
can be described uniformly in the variety of commutative diassociative loops
whose $2$-primary component is C and whose $3$-primary component is Moufang.
More importantly, the encompassing variety happens to be the variety of commutative
\emph{RIF loops}, i.e., inverse property loops satisfying either, and hence both,
of the following identities:
\[
(\RIF{1}) \qquad (xy\cdot z)\cdot xy = x\cdot y(zx\cdot y), \qquad\qquad
(\RIF{2}) \qquad xy\cdot (z\cdot xy) = (x\cdot yz)x\cdot y \,.
\]
These loops were defined for the first time in \cite{KKP}.

To understand the structure of commutative RIF loops requires the
study of \emph{Moufang elements}. These are traditionally defined
(for well-motivated reasons) to be those elements $x$ satisfying either
of the top two equations of (\Mfg) for every $y,z$. However, they
could certainly be defined in other natural and non-equivalent ways,
by fixing any variable in any one of the equations in (\Mfg), and
assuming that the other two variables in that equation are universally
quantified.

We analyze the situation in \S\secref{elements}, which we hope will eventually
lead to a deeper understanding of Moufang elements. We could not
resist the temptation and proved somewhat more than is needed for \S\secref{rif},
but the topic remains rife with open problems, some of which we state explicitly.

The main results of this paper can be found in \S\secref{rif}, where we
describe the structure of commutative RIF loops and give the main decomposition
theorem.

Finally, it is well-known that commutative Moufang loops
are closely related to totally symmetric quasigroups, and commutative C loops
to Steiner triple systems. As an application of our results, we conclude the paper
in \S\secref{quasigroups} by showing how commutative RIF loops are related to a
certain class of quasigroups, recovering the C and Moufang situations as
special cases.

Our investigations were aided by the automated theorem prover Prover9
\cite{Prover9}, the finite model builder Mace4 \cite{Mace4}, and the
LOOPS package \cite{LOOPS} for GAP \cite{GAP}.

\section{The general decomposition}
\seclabel{general}

A subloop $N$ of a loop $Q$ is \emph{normal}, denoted $N\unlhd Q$, if it is
a kernel of some loop homomorphism with domain $Q$. When $S$ is a subset of
$Q$, we let $\langle S\rangle$ denote the subloop of $Q$ generated by $S$.

Let $\setof{Q_i}{i\in I}$ be a collection of subloops of a loop $Q$.
Then $Q$ is the \emph{(internal) direct product} of $\setof{Q_i}{i\in I}$
if
\begin{enumerate}
\item[(i)]\quad $Q_i\unlhd Q$ for every $i$,
\item[(ii)]\quad $Q_i\cap \sbl{Q_j\,|\,j\in I,\,j\ne i} = 1$,
\item[(iii)]\quad $Q = \sbl{Q_i\,|\,i\in I}$.
\end{enumerate}
If the index set $I$ is finite, the internal direct product $Q$ of
$\setof{Q_i}{i\in I}$ is isomorphic to the \emph{external direct product}
$\prod_{i\in I} Q_i$, where multiplication is performed componentwise
(\cite[Lemma IV 5.1]{Bruck}).

For a power-associative loop $Q$ and a positive integer $k$, let $Q_{[k]}$
denote the set of all torsion elements $x\in Q$ such that $|x|$ divides $k$.

\begin{lemma}
\lemlabel{dec1}
Let $Q$ be a commutative diassociative loop.
\begin{enumerate}
\item[(i)]\quad For each $n\geq 0$, the mapping
$Q\to Q; x\mapsto x^n$ is a homomorphism
with kernel $Q_{[n]}$.
\item[(ii)]\quad For any torsion elements $x_1, \ldots, x_k\in Q$,
$|x_1\cdots x_k|$ is a divisor of $\lcm \{|x_1|, \ldots, |x_k|\}$,
no matter how $x_1\cdots x_k$ is parenthesized.
\end{enumerate}
\end{lemma}

\begin{proof}
We have $(xy)^n = x^n y^n$ immediately from commutativity and diassociativity,
and so (i) follows. If $x_1, \ldots, x_k$ are torsion elements, let $n = \lcm
\{|x_1|,\ldots,|x_k|\}$. Then $(x_1\cdots x_k)^n = x_1^n\cdots x_k^n$, where
the two products are parenthesized in analogous way. Since $x_j^n = 1$ for each
$j$, we have (ii).
\end{proof}

For each $x$ in a loop $Q$, the \emph{left translation} $L_x$ and
the \emph{right translation} $R_x$ are permutations of $Q$ defined,
respectively, by $L_x y := xy$ and $R_x y := yx$ for all $y\in Q$.
The \emph{inner mapping group} $\Inn(Q)$ of a loop $Q$ is the
stabilizer of the neutral element $1$ in the group
generated by all left and right translations. $\Inn(Q)$ is
generated by all permutations of the forms $R_x\inv L_x$,
$L_{xy}\inv L_x L_y$ and $R_{yx}\inv R_x R_y$ \cite{Bruck}.

Recall that a subloop $P\leq Q$ is normal in $Q$ if and only if $\varphi P
\subseteq P$ for all $\varphi \in \Inn(Q)$, that is, if and only if $P$ is
invariant under the action of $\Inn(Q)$. With this characterization of
normality, the following is obvious.

\begin{lemma}
\lemlabel{tower}
Let $\{P_i\}_{i=1}^{\infty}$ be a sequence of normal subloops of
a loop $Q$ satisfying $P_i \leq P_{i+1}$ for each $i$. Then
$\bigcup_{i=1}^{\infty} P_i$ is a normal subloop.
\end{lemma}

\begin{lemma}
\lemlabel{dec2}
Let $Q$ be a commutative diassociative loop. Then
for each prime $p$, $\pc{Q}{p} \normal Q$.
\end{lemma}

\begin{proof}
By Lemma \lemref{dec1}, $Q_{[p^m]} \normal Q$ for each
$m\geq 0$. Since
$Q_{[p^m]} \leq Q_{[p^{m+1}]}$ for each $m$, and also
$\pc{Q}{p} = \bigcup_{m\geq 0} Q_{[p^m]}$,
we have $\pc{Q}{p}\normal Q$ by Lemma \lemref{tower}.
\end{proof}

\begin{lemma}
\lemlabel{factor}
Let $Q$ be a commutative diassociative loop. If $m,n$ are
relatively prime positive integers, then $Q_{[mn]} = Q_{[m]} Q_{[n]}$,
a direct product.
\end{lemma}

\begin{proof}
If $x\in Q_{[m]}$ and $y\in Q_{[n]}$, then $(xy)^{mn} = x^{mn} y^{mn} = 1$,
and so $Q_{[m]} Q_{[n]}\subseteq Q_{[mn]}$. Now fix $z\in Q_{[mn]}$ and choose
$r,s$ so that $mr + ns = 1$. Then $z = z^{ns} z^{mr}$. Since
$z^{ns}\in Q_{[m]}$ and $z^{mr}\in Q_{[n]}$, we have the other inclusion.
The product is direct because each $Q_{[j]}$ is normal (Lemma \lemref{dec1})
and $Q_{[m]}\cap Q_{[n]} = \{1\}$.
\end{proof}

\begin{theorem}[Bruck and Paige \cite{BruckPaige}]
\thmlabel{general}
A torsion, commutative diassociative loop is the direct
product of its $p$-primary components, that is, a direct
product of commutative diassociative $p$-loops.
\end{theorem}

\begin{proof}
Let $Q$ be a torsion, commutative diassociative loop. In view of Lemma
\lemref{dec2}, it remains to show that $Q=\sbl{\pc{Q}{p}\,|\,p\ \text{prime}}$,
and $\pc{Q}{p}\cap \sbl{\pc{Q}{q}\,|\,q\ne p,\,\text{$q$ a prime}} = 1$.

Fix $x\in Q$ with $x\neq 1$. Since $Q$ is torsion,
$x\in Q_{[n]}$ for some $n > 0$. By Lemma \lemref{factor} and
induction, $Q_{[n]} = Q_{[p_1^{a_1}]} \cdots Q_{[p_k^{a_k}]}$
(direct product)
where $n = p_1^{a_1}\cdots p_k^{a_k}$ for some distinct primes
$p_i$ and exponents $a_i > 0$. Since
$Q_{[p_i^{a_i}]} \subseteq \pc{Q}{p_i}$, we have
$x\in \pc{Q}{p_1}\cdots \pc{Q}{p_k}$. This shows
$Q=\sbl{\pc{Q}{p}\,|\,p\ \text{prime}}$.

Now assume that $x\in \pc{Q}{p}\cap \sbl{\pc{Q}{q}\,|\,q\ne p}$. Then
$x\in \sbl{\pc{Q}{q_1},\dots,\pc{Q}{q_k}}$ for some $q_i\ne p$. Since
all $\pc{Q}{q_i}$ are normal in $Q$ by Lemma \lemref{dec1}, we have $\langle
\pc{Q}{q_1},\dots,\pc{Q}{q_k}\rangle = \pc{Q}{q_1}\cdots \pc{Q}{q_k}$. Thus $x
= x_1\cdots x_k$, where $x_i\in \pc{Q}{q_i}$, $|x_i|=q_i^{a_i}$, and the
product $x_1\cdots x_k$ is parenthesized in some way. By Lemma \lemref{dec1},
$|x|$ is a divisor of $q_1^{a_1}\cdots q_k^{a_k}$. But $|x|$ is also a power
of $p$, so we conclude that $x=1$.
\end{proof}

The \emph{nucleus} and \emph{center} of a loop $Q$ are the sets
\begin{align*}
N(Q) &= \setof{a\in Q}{a\cdot xy = ax\cdot y, x\cdot ay = xa\cdot y,
x\cdot ya = xy\cdot a,\ \forall x,y\in Q}, \\
Z(Q) &= N(Q) \cap \setof{a\in Q}{ax = xa,\ \forall x\in Q} .
\end{align*}
The nucleus is a subloop of $Q$, but is not necessarily normal.
The center is a normal subloop of any loop.
In a commutative loop, the center and nucleus coincide.

Note that if $Q=\prod_i Q_i$ then $Z(Q) = \prod_i Z(Q_i)$.
It is now easy to see what happens if we impose the condition that $x^n$ is
central in torsion commutative diassociative loops.

\begin{theorem}
\thmlabel{central} Let $n>0$ be a fixed integer, and let $p_1^{a_1}\cdots
p_k^{a_k}$ be a prime factorization of $n$. Let $Q$ be a torsion commutative
diassociative loop with each $x^n\in Z(Q)$. Then $Q$ is a direct product of
commutative diassociative $p_i$-loops in which $p_i^{a_i}$th powers are central
with an abelian group in which each element has order prime to $n$.
\end{theorem}

\begin{proof}
Let $Q$ be a torsion commutative diassociative loop. Let $x\in \pc{Q}{p_i}$ and
$m = n/p_i^{a_i}$. Since $|x|=p_i^{b_i}$ for some $b_i$, $m$ and $|x|$ are
relatively prime, and so $x^m$ is a generator of $\langle x\rangle$. In
particular, $x=x^{rm}$ for some $r$. Thus $x^{p_i^{a_i}} = x^{rn}\in Z(Q)\cap
\pc{Q}{p_i}=Z(\pc{Q}{p_i})$, where the last equality holds because $Q$ is
a direct product of its $p$-primary components (Theorem \thmref{general}).

Conversely, let $Q$ be a direct product of an abelian group $G$ and
diassociative $p_i$-loops $Q_i$ in which $p_i^{a_i}$th powers are central. Then
clearly $x^n\in Z(G)\cap \bigcap_i Z(Q_i) = Z(Q)$.
\end{proof}

\section{Moufang elements}
\seclabel{elements}

There are various instances of diassociativity to which we will need
to make specific reference. The
\emph{inverse property} (IP) is defined by any two of the
following equations (which together imply the third):
\[
(\LIP) \qquad x\inv\cdot xy = y, \qquad\qquad
(\RIP) \qquad xy\cdot y\inv = x, \qquad\qquad
(\AAIP) \qquad (xy)\inv = y\inv x\inv\,.
\]
These are known, respectively, as the \emph{left inverse}, \emph{right
inverse}, and \emph{antiautomorphic inverse properties}.

\begin{remark} Not all loops have two-sided inverses. Given a loop $Q$ and
$x\in Q$, there are unique $x\linv , x\rinv\in Q$ such that $x\linv x = xx\rinv
= 1$. Then one can say that $Q$ has the inverse property if $x\linv \cdot xy =
yx\cdot x\rinv = y$ for all $x,y\in Q$. But these identities imply $x\linv =
x\rinv = x\inv$, so the inverse property can equivalently be stated as above.
Moreover, in the commutative case, which we deal with in \S\secref{rif}, we get
$x\linv = x\rinv = x\inv$ for free.
\end{remark}

We will also need the \emph{left alternative}, \emph{right alternative}, and
\emph{flexible} laws:
\[
(\LAlt) \qquad x\cdot xy = x^2 y, \qquad\qquad (\RAlt) \qquad xy\cdot y = xy^2,
\qquad\qquad (\Flex) \qquad x\cdot yx = xy\cdot x\,.
\]
Loops satisfying both (\LAlt) and (\RAlt) are called \emph{alternative}.

Moufang loops are RIF loops, but flexible C-loops are not necessarily RIF.
Both are included in a larger variety called \emph{ARIF loops}
(``Almost RIF''), which
are defined to be flexible loops satisfying either, and hence both, of
the identities
\[
(\ARIF{1}) \qquad x(yxy\cdot z) = xyx\cdot yz
\qquad\qquad
(\ARIF{2}) \qquad (z\cdot yxy)x = zy\cdot xyx
\]
These loops were introduced in \cite{KKP},
and the main result of that paper was the following.

\begin{proposition}
\prplabel{diassoc}
Every ARIF loop, and hence every RIF loop, is diassociative.
\end{proposition}

\noindent We will use Proposition \prpref{diassoc} freely throughout what follows.

Recall that an \emph{autotopism} of a loop $Q$ is a triple $(f,g,h)$ of
permutations of $Q$ satisfying $f(x)g(y)=h(xy)$ for all $x,y\in Q$. Observe:

\begin{proposition}
\prplabel{ip-aut}
Let $Q$ be an IP loop, let $J:Q\to Q;x\mapsto x\inv$ denote the inversion
mapping, and let $f,g,h$ be permutations of $Q$.
The following are equivalent:
\begin{enumerate}
\item[(i)] $(f,g,h)$ is an autotopism,
\item[(ii)] $(JfJ, h,g)$ is an autotopism,
\item[(iii)] $(h,JgJ,f)$ is an autotopism.
\end{enumerate}
\end{proposition}

We assume for the rest of this section that the \emph{flexible law holds}.
(We make this assumption to keep the situation manageable, although many of our
arguments would work without it, too.)

There are thus $3$ distinct Moufang identities (\Mfg), each with three variables.
We now consider elements defined by fixing a variable in a Moufang identity. In
anticipation of Lemma \lemref{mfg-equivs} below, we group the various
possibilities as follows:
\[
\begin{array}{rc}
(M_0) &
c\cdot xy\cdot c = cx\cdot yc, \qquad
c(x\cdot cy) = cxc\cdot y, \qquad
(yc\cdot x)c = y\cdot cxc, \\
(M_1) &
x(c\cdot xy) = xcx\cdot y, \qquad
(yx\cdot c)x = y\cdot xcx, \\
(M_2) &
xc\cdot yx = x\cdot cy\cdot x, \qquad
xyx\cdot c = x(y\cdot xc), \\
(M_3) &
xy\cdot cx = x\cdot yc\cdot x, \qquad
c\cdot xyx = (cx\cdot y)x.
\end{array}
\]
Each of these equations is assumed to be universally quantified in the
variables $x$ and $y$.

We can view these identities as nine possibly different definitions of
``Moufang elements.'' A natural question then is:

\begin{problem} What are all the implications among the nine definitions of
Moufang elements in the variety of flexible loops?
\end{problem}

Without additional assumptions, we are not able to establish a single
implication. However, in the IP case we have:

\begin{lemma}
\lemlabel{mfg-equivs}
For an element $c$ of a flexible IP loop $Q$,
\begin{enumerate}
\item[(i)] the equations ($M_0$) are equivalent,
\item[(ii)] the equations ($M_1$) are equivalent,
\item[(iii)] the equations ($M_2$) are equivalent,
\item[(iv)] the equations ($M_3$) are equivalent.
\end{enumerate}
\end{lemma}

\begin{proof}
For (i): In IP loops, we have $J L_x J = R_x\inv$. Now the three equations are
equivalent, respectively, to $(L_c,R_c,L_c R_c)$ being an autotopism, to $(R_c
L_c, L_c\inv,L_c)$ being an autotopism, and to $(R_c\inv, L_c R_c,R_c)$ being
an autotopism. The desired equivalence then follows from Proposition
\prpref{ip-aut} applied to $f = L_c$, $g = R_c$ and $h = L_c R_c$.

For (ii): If $xcx\cdot y = x(c\cdot xy)$ holds, then replace $y$ with $x\inv
(c\inv \cdot x\inv y\inv)= [(yx\cdot c)x]\inv$ to get $xcx\cdot [(yx\cdot
c)x]\inv = x(c\cdot x[x\inv (c\inv \cdot x\inv y\inv)]) = y\inv$. Thus $xcx =
y\inv\cdot (yx\cdot c)x$, and so $(yx\cdot c)x = y\cdot xcx$. The reverse
implication follows from the mirror of this argument.

For (iii): if $xc\cdot yx = x\cdot cy\cdot x$ holds, then $xc\cdot (c\inv \cdot
x\inv y\inv x\inv)x = x\cdot x\inv y\inv x\inv\cdot x = y\inv$, using the IP.
Thus $c\inv \cdot x\inv y\inv x\inv = (c\inv x\inv\cdot y\inv)x\inv$, and then
using (\AAIP) gives $xyx\cdot c = x(y\cdot xc)$. Conversely, if $xyx\cdot c =
x(y\cdot xc)$, then following the argument in reverse gives $xc\cdot (c\inv
\cdot x\inv y\inv x\inv)x = y\inv$. Replacing $y$ with $x\inv \cdot y\inv c\inv
\cdot x\inv$ and using (\AAIP) gives $xc\cdot yx = x\cdot cy\cdot x$.

Finally, the mirror of the proof of (iii) proves (iv).
\end{proof}

For a flexible IP loop $Q$, let $M_i(Q)$, $i=0,\ldots,3$ denote the sets of
elements satisfying, respectively, ($M_i$), $i=0,\ldots,3$. When the underlying
loop $Q$ is clear, as will usually be the case, we abbreviate $M_i = M_i(Q)$.

Elements of $M_0$ are known as \emph{Moufang elements}
(\cite{Bruck}, p. 113). This definition is motivated by isotopy considerations;
an element of an IP loop is contained in $M_0$ if and only if the loop
isotope defined by that element has the IP. See \cite{Bruck} for details.

\begin{lemma}
\lemlabel{inverse-equiv}
Let $Q$ be a flexible, IP loop. Then
\begin{enumerate}
\item[(i)] $M_0$ is a subloop,
\item[(ii)] $c\in M_1$ if and only if $c\inv\in M_1$,
\item[(iii)] $c\in M_2$ if and only if $c\inv\in M_3$.
\end{enumerate}
\end{lemma}

\begin{proof}
Part (i) is (\cite{Bruck}, Chap. VII, Lemma 2.2). The rest follows immediately
from (\AAIP).
\end{proof}

In a flexible IP loop which is not left alternative, the neutral element $1$
satisfies $1\in M_0$, but $1\not\in M_1$. The smallest order for which such a
loop exists is $7$ (this fact can be checked by computer with the help of any
library of small loops, for instance the one found in the GAP \cite{GAP}
package LOOPS \cite{LOOPS}, or with a model builder, such as \cite{Mace4}):
\begin{table}[htb]
\[
\begin{array}{c|ccccccc}
\cdot & 1 & 2 & 3 & 4 & 5 & 6 & 7 \\
\hline
1 & 1 & 2 & 3 & 4 & 5 & 6 & 7 \\
2 & 2 & 3 & 1 & 6 & 7 & 5 & 4 \\
3 & 3 & 1 & 2 & 7 & 6 & 4 & 5 \\
4 & 4 & 7 & 6 & 5 & 1 & 2 & 3 \\
5 & 5 & 6 & 7 & 1 & 4 & 3 & 2 \\
6 & 6 & 4 & 5 & 3 & 2 & 7 & 1 \\
7 & 7 & 5 & 4 & 2 & 3 & 1 & 6
\end{array}
\]
\end{table}
For instance, $5\cdot 1(5\cdot 6) = 5\cdot 3 = 7$, but
$(5\cdot 1\cdot 5)\cdot 6 = 4\cdot 6 = 2$.

\begin{lemma}
\lemlabel{alt-equiv}
In a flexible, alternative, IP loop,
$M_0 \subseteq M_1$\,.
\end{lemma}

\begin{proof}
If $c\in M_0$, then $c\inv\in M_0$ (since $M_0$ is a subloop),
and so for all $x,y$,
\begin{align*}
x\cdot yc\inv y &\textjeq{\LIP}
x\cdot [(c\cdot c\inv y)\cdot c\inv y]
\textjeq{\RAlt, \LIP}
c\inv\cdot c[x\cdot c(c\inv y)^2 ] \\
&\jeq{c\in M_0}
c\inv [ cxc\cdot (c\inv y)^2 ]
\textjeq{\RAlt}
c\inv [ (cxc\cdot c\inv y) \cdot c\inv y ] \\
& \jeq{c\in M_0}
c\inv [ (c\cdot x(c\cdot c\inv y)) \cdot c\inv y ]
\textjeq{\LIP}
c\inv [ (c\cdot xy) \cdot c\inv y ] \\
&\jeq{c\inv\in M_0}
c\inv (c\cdot xy) c\inv \cdot y
\textjeq{\LIP}
(xy\cdot c\inv )y \,.
\end{align*}
Thus $c\inv \in M_1$, and so
$c\in M_1$ (Lemma \lemref{inverse-equiv}).
\end{proof}

\begin{problem}
\prblabel{alt}
Does there exist a diassociative loop in which $M_0\neq M_1$?
A flexible, alternative, IP loop?
\end{problem}

\begin{lemma}
\lemlabel{two-hyp}
In a flexible, IP loop, $M_0\cap M_2 = M_0\cap M_3$.
\end{lemma}

\begin{proof}
Fix $c\in M_0\cap M_2$. We compute
\[
c[x\cdot yc\cdot x]c
\jeq{c\in M_0}
cx\cdot (yc\cdot x)c
\jeq{c\in M_0}
cx\cdot (y\cdot cxc)
\jeq{c\in M_2}
(cx\cdot y\cdot cx)c\,.
\]
Canceling $c$ on the right, and then multiplying on
the left by $c\inv$ and using (\LIP), we get
\[
x\cdot yc\cdot x = c\inv (cx\cdot y\cdot cx)
\jeq{c\inv \in M_3}
(c\inv\cdot cx)y\cdot cx
\textjeq{\LIP}
xy\cdot cx\,,
\]
where we have used Lemma \lemref{inverse-equiv} in
the second step. Thus $c\in M_3$.

Conversely, if $c\in M_0\cap M_3$, then
$c\inv\in M_0\cap M_2$ (Lemma \lemref{inverse-equiv}),
and so $c\inv\in M_3$ by the preceding paragraph.
Thus $c\in M_2$ (Lemma \lemref{inverse-equiv} again).
This completes the proof.
\end{proof}

\begin{problem}
\prblabel{m2m3}
Does there exist a diassociative loop in which
$M_2\neq M_3$? A flexible, alternative, IP loop?
A flexible IP loop?
\end{problem}

\begin{theorem}
\thmlabel{arif-mfg}
In an ARIF loop,
$M_2 = M_3 \subseteq M_0 = M_1$.
\end{theorem}

\begin{proof}
If $c\in M_1$, then
\[
c\cdot x(c\cdot y) \textjeq{\LIP}
c\cdot x(c\cdot [x\cdot x\inv y])
\jeq{c\in M_1}
c(xcx\cdot x\inv y) \textjeq{\ARIF{1}}
cxc\cdot [x\cdot x\inv y]
\textjeq{\LIP}
cxc\cdot y\,.
\]
Therefore $c\in M_0$.
We then have $M_0 = M_1$ by Lemma \lemref{alt-equiv}.

Now suppose $c\in M_3$. Then
\begin{align*}
y\inv xy\inv \cdot y(y\inv x\cdot c)y
&\jeq{c\in M_3} y\inv xy\inv \cdot [(y\cdot y\inv x)\cdot cy]
\textjeq{\LIP}
y\inv xy\inv \cdot [x\cdot cy] \\
&\textjeq{\ARIF{1}}
y\inv [xy\inv x\cdot cy]
\textjeq{\LIP, \RAlt}
y\inv [(y\cdot (y\inv x)^2)\cdot cy] \\
&\jeq{c\in M_3}
y\inv [y\cdot (y\inv x)^2 c\cdot y]
\textjeq{\LIP}
(y\inv x)^2 c\cdot y \,.
\end{align*}
Replacing $x$ with $yx$ and using (\LIP), we have
$xy\inv \cdot [y\cdot xc\cdot y] = x^2 c\cdot y$, and so
by (\RAlt),
\begin{align*}
(x\cdot xc) y\cdot xc &=
(xy\inv \cdot [y\cdot xc\cdot y])\cdot xc
\textjeq{\ARIF{2}}
x(xc\cdot y\cdot xc) \\
&\textjeq{\RIP}
(xc\cdot c\inv)[xc\cdot y\cdot xc]
\jeq{c\inv \in M_2}
xc \cdot c\inv (xc\cdot y)\cdot xc \,,
\end{align*}
where we use Lemma \lemref{inverse-equiv} in the last step.
Canceling $xc$ on the right and replacing
$x$ with $xc\inv$ and using (\RIP), we get
$xc\inv x\cdot y = x(c\inv\cdot xy)$. Thus
$c\inv\in M_1$, and so $c\in M_1$ by Lemma \lemref{inverse-equiv}.
This establishes $M_3\subseteq M_1 = M_0$.
By Lemma \lemref{inverse-equiv}, we thus also have
$M_2\subseteq M_1 = M_0$.  By Lemma \lemref{two-hyp},
$M_2 = M_3$. This completes the proof.
\end{proof}

\begin{problem}
\prblabel{arif}
Is there an ARIF loop in which $M_2 \neq M_1$?
\end{problem}

Problem \prbref{arif} has a negative answer for
the two major subvarieties of the ARIF variety,
namely RIF loops (Theorem \thmref{rif-mfg}) and
flexible C-loops (Corollary \corref{flexc-mfg}).

\begin{theorem}
\thmlabel{rif-mfg}
In a RIF loop, $M_i = M_j$ for all $i,j\in \{0,1,2,3\}$.
\end{theorem}

\begin{proof}
In view of Theorem \thmref{arif-mfg}, it
is enough to show that $M_0\subseteq M_2$.

In the RIF identity (\RIF{1}), replace $x$ with $xy\inv$,
use (\RIP), and then replace $y$ with $y\inv$ to get
$xzx = xy\cdot [y\inv (z\cdot xy) y\inv]$. Now assume $c\in M_0$,
and set $z = cu$ and $y = c$ to obtain
\begin{align*}
x\cdot cu\cdot x &= xc\cdot (c\inv (cu \cdot xc) c\inv)
\jeq{c\in M_0}
xc\cdot (c\inv (c\cdot ux\cdot c) c\inv)
\textjeq{\LIP, \RIP}
xc\cdot ux\,.
\end{align*}
Thus $c\in M_2$.
\end{proof}

An element $c$ of a loop $Q$ is a \emph{C element} if it satisfies
\[
x(c\cdot cy) = (xc\cdot c)y
\]
for all $x,y\in Q$. Note that C elements satisfy $c\cdot cx = c^2 x$ and
$xc\cdot c = xc^2$ for all $x\in Q$, which we use below without reference.

Let $C_0 = C_0(Q)$ denote the set of all C elements of $Q$. Chein \cite{Chein}
showed the following:

\begin{proposition}
\prplabel{c-el}
In an IP loop $Q$, $c\in C_0(Q)$ if and
only if $c^2\in N(Q)$.
\end{proposition}

\begin{lemma}
\lemlabel{c-mfg}
In a flexible IP loop,
$C_0\cap M_0 \subseteq C_0\cap M_2 = C_0\cap M_3$.
\end{lemma}

\begin{proof}
If $c\in C_0\cap M_2$, then $c\inv \in M_3$
(Lemma \lemref{inverse-equiv}), and so
\[
(cx\cdot y)x \textjeq{\LIP}
(c^2\cdot c\inv x)y\cdot x \jeq{c^2\in N}
c^2\cdot (c\inv x\cdot y)x
\jeq{c\inv \in M_3}
c^2 \cdot (c\inv \cdot xyx)
\textjeq{\LIP}
c\cdot xyx \,.
\]
Thus $c\in M_3$. Therefore
$C_0\cap M_2\subseteq C_0\cap M_3$ and
Lemma \lemref{inverse-equiv} gives the reverse inclusion.

Now suppose $c\in C_0\cap M_0$. Then
\[
cx\cdot yc\cdot cx \jeq{c\in M_0}
(c\cdot xy\cdot c)\cdot cx \jeq{c\in M_0}
c (xy \cdot c^2 x) \jeq{c^2\in N}
c (x\cdot yc^2\cdot x)\,.
\]
Replace $y$ with $yc\inv$, use (\RIP), and multiply
on the right by $c$:
\[
(cx\cdot y\cdot cx)c = c(x\cdot yc\cdot x)c
\jeq{c\in M_0}
cx\cdot (yc\cdot x)c \jeq{c\in M_0}
cx\cdot y(cx\cdot c)\,.
\]
Replace $x$ with $c\inv x$ and use (\LIP) to get
$xyx\cdot c = x(y\cdot xc)$, that is, $c\in M_2$.
\end{proof}

\begin{theorem}
\thmlabel{flexc-mfg}
In an ARIF loop, $C_0\cap M_i = C_0\cap M_j$ for all
$i,j\in \{0,1,2,3\}$.
\end{theorem}

\begin{proof}
This follows immediately from
Theorem \thmref{arif-mfg} and
Lemma \lemref{c-mfg}.
\end{proof}

\begin{corollary}
\corlabel{flexc-mfg}
In a flexible C-loop, $M_i = M_j$ for all
$i,j\in \{0,1,2,3\}$.
\end{corollary}

\section{Commutative RIF loops}
\seclabel{rif}

We begin with some characterizations of the variety of
commutative RIF loops.

\begin{lemma}
\lemlabel{char1}
A loop $Q$ is a commutative RIF loop if and only if it is an IP
loop satisfying the identity
\[
x (y^2\cdot xz) = (xy)^2 z\,.  \tag{\CRIF}
\]
for all $x,y,z\in Q$.
\end{lemma}

\begin{proof}
In a commutative, alternative loop, we have
\[
(xy\cdot z)\cdot xy = xy\cdot (xy\cdot z) = (xy)^2 z
\qquad \text{and}\qquad
x\cdot y(zx\cdot y) = x\cdot y(y\cdot xz) = x(y^2\cdot xz)\,.
\]
In RIF loops, the left hand sides are equal, and since such loops
are diassociative, it follows that commutative RIF
loops satisfy (\CRIF). To complete the proof, it is enough to show that
an IP loop satisfying (\CRIF) is alternative and commutative. Taking
$y = 1$ in (\CRIF), we get $x\cdot xz = x^2 z$ which is (\LAlt).
By (\AAIP), any identity in an IP loop is equivalent to its
mirror, so we also have $zx\cdot x = z x^2$, that is, (\RAlt).
Taking $z = 1$ in (\CRIF) gives $x\cdot y^2 x = (xy)^2$, which is equivalent
to $y^2 x = x\inv (xy)^2$. Applying (\LAlt) and (\RAlt),
we have $y\cdot yx = (x\inv\cdot xy)\cdot xy = y\cdot xy$. Canceling,
it follows that $Q$ is commutative.
\end{proof}

The identity (\CRIF) has appeared in the literature before
in other contexts. It plays a role in the theory of, for
instance, Bruck loops \cite{Kiechle}.

The following is evidence of the naturality of the variety of
commutative RIF loops. Among other things, it shows that
passing from RIF to ARIF adds no generality in the commutative case.

\begin{theorem}
\thmlabel{crif-char}
For a commutative loop $Q$, the following are equivalent.
\begin{enumerate}
\item[(i)] $Q$ is a RIF loop,
\item[(ii)] $Q$ is an ARIF loop.
\item[(iii)] $Q$ is an alternative, IP loop
with each $x^2\in M_0(Q)$,
\item[(iv)] $Q$ satisfies $(\CRIF)$,
\end{enumerate}
\end{theorem}

\begin{proof}
(i)$\implies$(ii) holds even in the noncommutative case \cite{KKP}.

For (ii)$\implies$(iii): By Proposition \prpref{diassoc}, $Q$ is diassociative,
and we freely use this and commutativity in the following calculation:
\begin{align*}
(z^2\cdot x^2 y)x\inv &= L_{z^2\cdot x^2 y} L_{x^2 y}(y\inv x^{-3})
\textjeq{ARIF{1}} L_z L_{z(x^2 y)^2} (x^3 y)\inv =
L_z L_{x^3 y}\inv L_{(x^2 y)^2} z \\
&= L_z L_{x^3 y}\inv L_{(x^3 y)^2 x^{-2}} L_{x^{-2}} (x^2 z)
\textjeq{ARIF{1}} L_z L_{x^3 y}\inv L_{x^3 y} L_{x^3 y\cdot x^{-4}} (x^2 z) \\
&= L_z L_{x\inv y} (x^2 z) = L_z L_{x^2 z} (x\inv y)
\textjeq{ARIF{1}} L_{z^2 x} L_x (x\inv y) = z^2 x\cdot y\,.
\end{align*}
Thus $z^2\cdot x^2 y = (z^2 x\cdot y)x$, and so each $z^2\in M_3$.
By Theorem \thmref{arif-mfg}, each $z^2 \in M_0$.

For (iii)$\implies$(iv): If each $y^2\in M_0$, then by Lemma
\lemref{alt-equiv}, each $y^2\in M_1$, and so $x (y^2\cdot xz) = x y^2 x\cdot
z$ for all $x,y,z\in Q$. By commutativity and the alternative laws, $xy^2 x =
x(y\cdot xy)$. Now $x\inv (xy)^2 \textjeq{\RAlt} (x\inv \cdot xy)\cdot xy
\textjeq{\LIP} y\cdot xy$, and so by (\LIP), $(xy)^2 = x(y\cdot xy)$. Thus
$x(y^2\cdot xz) = (xy)^2 z$ for all $x,y,z\in Q$, that is, (\CRIF) holds.

For (iv)$\implies$(i): take $y=1$ in (\CRIF) to get (\LAlt), and by
commutativity, (\RAlt). Also,
\[
xy\cdot y \textjeq{\RAlt} xy^2 = x(y^2\cdot xx\inv)
\textjeq{\CRIF} (xy)^2 x\inv
\textjeq{\LAlt} xy\cdot (xy\cdot x\inv)\,.
\]
Canceling and using commutativity, we obtain $y = x\inv\cdot xy$, and so the IP
holds. By Lemma \lemref{char1}, (i) holds.
\end{proof}

The following is well-known, and holds in more generality than we give here.

\begin{lemma}
\lemlabel{comm-mf} Let $Q$ be a commutative, IP loop. Then for every $x\in
M_0(Q)$, $x^3 \in Z(Q)$.
\end{lemma}

\begin{proof}
By (\cite{Bruck}, Chap. VII, Lemma 2.2), in an IP loop, for each $x\in M_0(Q)$,
the inner mapping $R_x\inv L_x$ is a pseudoautomorphism with companion
$x^{-3}$, that is, $x^{-3} \cdot (x \cdot yz) x\inv = [x^{-3} \cdot (x \cdot y)
x\inv]\cdot (x \cdot z) x\inv$, for all $y,z\in Q$. In the commutative case,
this reduces to $x^{-3} \cdot yz = x^{-3} y \cdot z$, that is, $x^{-3}\in
Z(Q)$.
\end{proof}

\begin{corollary}
\corlabel{exp6}
Let $Q$ be a commutative RIF loop. Then $Q/Z(Q)$ has exponent $6$.
\end{corollary}

\begin{proof}
By Theorem \thmref{crif-char}, every square is a Moufang element.
Then by Lemma \lemref{comm-mf}, every sixth power is central.
\end{proof}

Recall that a \emph{Steiner loop} is an IP loop of exponent $2$, or
equivalently, a C loop of exponent $2$ \cite{PhillipsVojtechovsky}.
Such loops are commutative.

\begin{theorem}
\thmlabel{mf-sq} Let $Q$ be a commutative RIF loop. Then:
\begin{enumerate}
    \item[(i)] For each $x\in Q$, $x^2\in  M_0(Q)$.
    \item[(ii)] $M_0(Q)$ is a normal subloop of $Q$.
    \item[(iii)] $Q/M_0(Q)$ is a C loop of exponent $2$, i.e., a Steiner
    loop.
\end{enumerate}
\end{theorem}
\begin{proof}
Part (i) is Theorem \thmref{crif-char}(iii). For (ii): The set of Moufang
elements is a subloop of any IP loop, so for (ii), only the normality requires
a proof. Fix $b,c\in Q$, $a\in M_0(Q)$, and set $d = L_{bc}\inv L_b L_c a =
(bc)\inv (b\cdot ca)$. We wish to show that $d\in M_0(Q)$. First, we compute
\[
b\cdot ca^2 = b[a^2 c^2 \cdot c\inv] = b[(ac)^2 \cdot b(bc)\inv]
= (b\cdot ac)^2 (bc)\inv = (bc)\inv (bc\cdot d)^2 = d(bc\cdot d) = bc\cdot d^2\,,
\]
where we have used (\CRIF) in the third equality, and commutativity and
diassociativity throughout. Thus $L_{bc}\inv L_b L_c (a^2) = d^2$. Now in RIF
loops, inner mappings preserve inverses \cite{KKP}, and so $b\cdot ca^{-2} =
bc\cdot d^{-2}$. Thus using $a^{-3}\in Z(Q)$ (Lemma \lemref{comm-mf}), we have
\[
a^{-3} d = (bc)\inv\cdot a^{-3}(b\cdot ca) =
(bc)\inv (b\cdot ca^{-2}) = (bc)\inv (bc\cdot d^{-2}) = d^{-2} \,.
\]
Therefore,
$d^3 = a^3 \in Z(Q)$. On the other hand, $d^2\in M_0(Q)$, and since
$Z(Q)\subseteq M_0(Q)$, we have $d = d^3 d^{-2}\in M_0(Q)$. This
completes the proof of normality.

Part (iii) then follows from (ii) and Theorem \thmref{crif-char}.
\end{proof}

Next we turn to C elements.
Although it is a bit of an aside to the rest of the development,
we mention the following in passing.

\begin{theorem}
\thmlabel{z=cm}
Let $Q$ be a commutative IP loop. Then $M_0\cap C_0 = Z(Q)$.
\end{theorem}

\begin{proof}
If $a\in M_0\cap C_0$, then $a^3\in Z$ (Lemma \lemref{comm-mf})
and $a^2\in Z$ (Proposition \prpref{c-el}), and so
$a = a^3 a^{-2}\in Z$. The other inclusion is clear.
\end{proof}

For commutative RIF loops, the subset of C elements is
well-structured.

\begin{theorem}
\thmlabel{c-normal} Let $Q$ be a commutative RIF loop. Then:
\begin{enumerate}
    \item[(i)] For each $x\in Q$, $x^3 \in C_0(Q)$.
    \item[(ii)] $C_0(Q)$ is a normal subloop of $Q$.
    \item[(iii)] $Q/C_0(Q)$ is a Moufang loop of exponent $3$.
\end{enumerate}
\end{theorem}

\begin{proof}
Part (i) follows from Proposition \prpref{c-el} and Corollary \corref{exp6}.

Now for $a,b\in C_0(Q)$, $(ab)^2 = a^2 b^2$ by diassociativity, and so by
Proposition \prpref{c-el}, $ab\in C_0(Q)$. In addition, $a\inv$ is clearly in
$C_0(Q)$, and so $C_0(Q)$ is a subloop. To show normality, fix $a\in C_0(Q)$,
$b,c\in Q$, and set $d = L_{bc}\inv L_b L_c a = (bc)\inv (b\cdot ca)$. We wish
to show $d\in C_0(Q)$. In RIF loops, inner mappings preserve inverses
\cite{KKP}, and so $b\cdot ca\inv = bc\cdot d\inv$. Using this and $a^{-2}\in
Z(Q)$, we compute
\[
a^{-2} d = (bc)\inv \cdot a^{-2}(b\cdot ca) = (bc)\inv (b\cdot ca\inv) =
(bc)\inv\cdot (bc\cdot d\inv) = d\inv \,.
\]
Thus $d^2 = a^2\in Z(Q)$. Since $d^3\in C_0(Q)$
and $Z(Q)\subseteq C_0(Q)$, we have $d = d^3 d^{-2}\in C_0(Q)$. This
completes the proof of (ii).

Finally, $Q/C_0(Q)$ has exponent $3$ by (i), and so by
Theorem \thmref{mf-sq}, every element of $Q/C_0(Q)$ is Moufang. This
proves (iii).
\end{proof}

Finally, we have our decomposition theorem in the torsion case.

\begin{theorem}
\thmlabel{rif-decomp}
Let $Q$ be a torsion, commutative RIF loop. Then $Q$ is the
direct product
of a C\ $2$-loop, a Moufang $3$-loop, and an abelian group
in which each element has order prime to $6$.
\end{theorem}

\begin{proof}
By Corollary \corref{exp6}, every sixth power is central. By
Theorem \thmref{central}, $Q$ is the direct product of a $2$-loop,
a $3$-loop, and an abelian group in which each element has order
prime to $6$. Since every cube is a C element
(Theorem \thmref{c-normal}), the $2$-primary component is a C loop.
Since every square is a Moufang element (Theorem \thmref{mf-sq}),
the $3$-primary component is Moufang.
\end{proof}

\section{Quasigroups associated to commutative RIF loops}
\seclabel{quasigroups}

Throughout this section, we will use multiplicative notation
for quasigroups, and additive notation for loops.
In particular, $0$ is the neutral element,
$-x$ is the inverse of $x$, and $x-y$ stands for $x+(-y)$ in loops.

A quasigroup $(Q,\cdot)$ is \emph{totally symmetric} if it is commutative and
satisfies the identity
\begin{displaymath}
    x\cdot xy = y \tag{TS}
\end{displaymath}
for every $x,y\in Q$. An element $0\in Q$ is an \emph{idempotent} if
$0^2 = 0$. Let \TS denote the category of totally symmetric quasigroups
with a distinguished idempotent element (uniformly denoted by $0$)
preserved by morphisms.

A loop $(Q,+)$ with two-sided inverses has the \emph{weak inverse property}
if it satisfies the identity
\begin{displaymath}
    x-(y+x) = -y\tag{WIP}
\end{displaymath}
for every $x,y\in Q$. Let \CWIP denote the category of commutative WIP loops.

Given a commutative quasigroup $(Q,\cdot)$ with an idempotent
$0\in Q$, define $\Loop(Q,\cdot) = (Q,+)$ by
\begin{displaymath}
    x + y = 0x\cdot 0y.
\end{displaymath}
Conversely, given a loop $(Q,+)$ with neutral element $0$, define
$\Quasi(Q,+)=(Q,\cdot)$ by
\begin{displaymath}
    x\cdot y = -x-y.
\end{displaymath}
It is then easy to show:

\begin{proposition} $\Quasi$ is a functor \TS$\to$\CWIP, and $\Loop$ is a functor
\CWIP$\to$\TS. Moreover, $\Loop\Quasi$ is identical on \TS, and $\Quasi\Loop$ is
identical on \CWIP, so that the categories \TS, \CWIP are equivalent.
\end{proposition}

The equivalence of \TS and \CWIP takes on a particularly nice form when
restricted to certain subcategories.

In a quasigroup $Q$ let $\Idem(Q)$ denote the set of all idempotents
of $Q$. In general, $\Idem(Q)$ need not be a subquasigroup of $Q$.
A quasigroup is said to be \emph{idempotent} if $Q = \Idem(Q)$.

A quasigroup $(Q,\cdot)$ is \emph{distributive} if it satisfies
\begin{displaymath}
    x(yz) = xy\cdot xz,\quad\quad (xy)z = xz\cdot yz\tag{D}
\end{displaymath}
for every $x$, $y$, $z\in Q$. Distributive quasigroups are idempotent.
The following result is due to Bruck \cite{Bruck44} (see also
\cite{Pflugfelder}, Thm. V.2.16).

\begin{proposition}
\prplabel{tsdist-cml3}
Let $(Q,+)$ be a commutative Moufang loop of exponent $3$. Then $\Quasi(Q,+)$
is a totally symmetric, distributive quasigroup. Conversely, let $(Q,\cdot)$ be
a totally symmetric, distributive quasigroup with a distinguished idempotent
$0$. Then $\Loop(Q,\cdot)$ is a commutative Moufang loop of exponent $3$.
\end{proposition}

A quasigroup is said to be \emph{unipotent} if $x^2=y^2$ for every $x,y$.

\begin{proposition}
\prplabel{steiner-steiner}
Let $(Q,+)$ be a $C$ loop of exponent $2$, i.e., a Steiner loop. Then
$\Quasi(Q,+)$ is an unipotent, totally symmetric quasigroup. Conversely, let
$(Q,\cdot)$ be a unipotent, totally symmetric quasigroup. Then $\Loop(Q,\cdot)$
is a Steiner loop.
\end{proposition}

Note that in a unipotent quasigroup there is a unique idempotent, namely
$0=x^2=y^2$. In a unipotent, totally symmetric quasigroup, it is easy to
see that the unique idempotent is a neutral element.
Thus the equivalence of Proposition \prpref{steiner-steiner} is
purely syntactical, since a unipotent, totally symmetric quasigroup \emph{is} a
Steiner loop. Put another way, the intersection of \TS and \CWIP
is precisely the variety of Steiner loops with neutral $0$, and
each of the functors $\Quasi$ and $\Loop$ is identical on that intersection.

Our task is to generalize simultaneously Propositions
\prpref{tsdist-cml3} and \prpref{steiner-steiner} by
finding the quasigroup counterpart of commutative RIF loops of
exponent $6$ under the functor $\Quasi$.

We introduce the following quasigroup axioms,
\begin{align*}
    x^2 x^2 &= x^2,  \tag{Q1} \\
    x (y^2\cdot xz) &= (xy)^2 z, \tag{Q2}
\end{align*}
noting that (Q2) is just another name for (\CRIF), this time in
quasigroups.

\begin{lemma}
\lemlabel{idempotent} A totally symmetric quasigroup satisfying \emph{(Q2)} is
distributive if and only if it is idempotent.
\end{lemma}

\begin{proof}
Only the sufficiency requires a proof. In the idempotent case,
(Q2) is equivalent to $x(y\cdot xz) = xy\cdot z$. Replacing $z$
with $xz$ and applying (TS), we obtain (D).
\end{proof}

For $0\in \Idem(Q)$, let $\Uni_0(Q) = \setof{x\in Q}{x^2 = 0}$.

\begin{lemma}
\lemlabel{squarehom}
Let $Q$ be a totally symmetric quasigroup satisfying \emph{(Q1)},
\emph{(Q2)}. Then
\begin{enumerate}
\item[(i)] the squaring mapping $Q\to Q; x\mapsto x^2$ is an
endomorphism of $Q$ with image $\Idem(Q)$,
\item[(ii)] $\Idem(Q)$ is a distributive subquasigroup,
\item[(iii)] for each $0\in \Idem(Q)$, $\Uni_0(Q)$ is a unipotent subquasigroup,
that is, a Steiner loop.
\end{enumerate}
\end{lemma}

\begin{proof}
Set $z = x$ in (Q2) and cancel $x$ on both sides to obtain $x^2 y^2 = (xy)^2$.
Thus squaring is an endomorphism. The image is a subset of $\Idem(Q)$ by (Q1),
and since every idempotent is trivially a square, the image
coincides with $\Idem(Q)$. This establishes (i). Homomorphic images of
quasigroups are subquasigroups, so (ii) follows from
Lemma \lemref{idempotent}. Finally, (iii) follows from (i).
\end{proof}

\begin{theorem}
\thmlabel{quasi-crif}
Let $(Q,+)$ be a commutative RIF loop of exponent $6$.
Then $(Q,\cdot) = \Quasi(Q,+)$ is a totally symmetric quasigroup satisfying
\emph{(Q1)}, \emph{(Q2)}.

Conversely, let $(Q,\cdot)$ be a totally symmetric quasigroup satisfying
\emph{(Q1)}, \emph{(Q2)}, and let $0\in Q$ be an idempotent. Then
$(Q,+) = \Loop(Q,\cdot)$ is a commutative RIF loop of exponent $6$.
\end{theorem}

\begin{proof}
Let $(Q,+)$ be a commutative RIF loop of exponent $6$, and let $(Q,\cdot) =
\Quasi(Q,+)$. Note that $x^2 = -2x$. Using diassociativity and the fact
that $Q$ has exponent $6$, we compute
\[
x^2\cdot x^2 = -(-2x) -(-2x) = 4x = -2x = x^2
\]
for all $x,y\in Q$. Thus (Q1) holds. Next
\begin{align*}
x(y^2\cdot xz) &= -x -(2y - (-x -z))
= -x + [-2y + (-x-z)] \\
&\textjeq{(\CRIF)} 2(-x-y) - z = (xy)^2 z
\end{align*}
for all $x,y,z\in Q$, and so (Q2) holds.

Now let $(Q,\cdot)$ be a totally symmetric quasigroup satisfying (Q1), (Q2),
and let $0\in Q$ be an idempotent. (Idempotents exist by (Q1).) Let
$(Q,+)=\Loop(Q,\cdot)$. First note that $2x = (0x)^2 = 0 x^2$ for all $x$ by
Lemma \lemref{squarehom}(i). We use this in the following calculations.
We verify (\CRIF) as follows:
\begin{align*}
x + [2y + (x + z)] &=  0x\cdot 0[(0\cdot 0y^2) \cdot 0(0x\cdot 0z)]
\textjeq{(TS)} 0x\cdot 0[y^2 \cdot 0(0x\cdot 0z)] \\
&\textjeq{(Q2)} 0x\cdot (0y)^2 (0x\cdot 0z) \textjeq{(Q2)} (0x\cdot 0y)^2 \cdot
0z
\textjeq{(TS)} [0\cdot 0(0x\cdot 0y)^2] \cdot 0z \\
&= 2(x + y) + z\,,
\end{align*}
where we used  in the second step. By Theorem \thmref{crif-char}, $(Q,+)$ is a
commutative RIF loop. We have
\[
6x = 2x + (2x + 2x) = (0\cdot 0x^2)\cdot 0[(0\cdot 0x^2)\cdot (0\cdot 0x^2)]
\textjeq{(TS)} x^2 \cdot 0 [x^2 x^2] \textjeq{(Q1)} x^2 \cdot 0 x^2
\textjeq{(TS)} 0,
\]
and this completes the proof.
\end{proof}

\begin{theorem}
\thmlabel{quasi-decomp}
Let $Q$ be a totally symmetric quasigroup satisfying \emph{(Q1)}, \emph{(Q2)}.
Then for each $0\in \Idem(Q)$, $Q$ is a direct product of $\Idem(Q)$ and
$\Uni_0(Q)$. Thus every totally symmetric quasigroup satisfying
\emph{(Q1)}, \emph{(Q2)} is a direct product of a distributive
subquasigroup and a Steiner loop.
\end{theorem}

\begin{proof}
Let $(Q,+)=\Loop(Q,\cdot)$ be the associated commutative RIF loop of exponent
$6$ (Theorem \thmref{quasi-crif}). By Theorem \thmref{rif-decomp}, $(Q,+)$ is a
direct product of a Moufang subloop $(Q_1,+)$ of exponent $3$ and a C subloop
$(Q_2,+)$ of exponent $2$, that is, a Steiner loop.

The subquasigroup $(Q_1,\cdot) = \Quasi(Q_1,+)$ of $(Q,\cdot)$ is distributive
(Proposition \prpref{tsdist-cml3}) and hence, idempotent. Thus $Q_1 \subseteq
\Idem(Q)$. On the other hand, if $c\in \Idem(Q)$, then
\begin{displaymath}
    c + (c + c) = 0c\cdot 0(0c\cdot 0c) = 0c\cdot 0(0^2 c^2) \textjeq{(TS)}
    0c\cdot c \textjeq{(TS)} 0,
\end{displaymath}
and so $c\in Q_1$. Therefore $Q_1 = \Idem(Q)$.

Next, the subquasigroup $(Q_2,\cdot) = \Quasi(Q_2,+)$ is just
$(Q_2,+)$ itself in different notation. In particular,
$Q_2\subseteq \Uni_0(Q)$. On the other hand, if $x\in \Uni_0(Q)$,
then $x^2 = 0$, and so $x + x = 0$, that is, $x\in Q_2$. Therefore
$Q_2 = \Uni_0(Q)$.

Finally, noting that the functor $\Quasi$ sends a direct product of
commutative diassociative loops to a direct product of quasigroups,
we have the desired result.
\end{proof}

\begin{remark} \emph{Steiner quasigroups} are defined as idempotent, totally
symmetric quasigroups. There is a one-to-one correspondence between Steiner
quasigroups of order $n$ and Steiner loops of order $n+1$. (Given a Steiner
quasigroup, introduce a new element $1$, leave $x\cdot y$ intact for $x\ne y$,
and set $x^2=1$. Conversely, given a Steiner loop with neutral element $1$,
remove $1$, leave $x\cdot y$ intact for $x\ne y$, and set $x^2=x$.) Moreover,
it is well-known that Steiner quasigroups are in one-to-one correspondence to
Steiner triple systems. Are there interesting combinatorial structures
associated to commutative RIF loops?
\end{remark}

\bibliographystyle{plain}

\begin{thebibliography}{99}

\bibitem{Bruck44}
R.~H.~Bruck,
Some results in the theory of quasigroups,
\textit{Trans. Amer. Math. Soc.} \textbf{55} (1944), 19--52.

\bibitem{Bruck}
R.~H.~Bruck,
\textit{A Survey of Binary Systems},
Springer-Verlag, 1971.

\bibitem{BruckPaige}
R.~H.~Bruck and L.~J.~Paige,
Loops whose inner mappings are automorphisms,
\textit{Ann. of Math.} (2) \textbf{63} (1956), 308--323.

\bibitem{Chein}
O.~Chein,
A short note on supernuclear (central) elements of inverse property loops.
\textit{Arch. Math. (Basel)} \textbf{33} (1979/80), 131--132.

\bibitem{GAP}
The GAP Group,
\textit{GAP -- Groups, Algorithms, and Programming, Version 4.4.9}; 2006.
\texttt{http://www.gap-system.org}

\bibitem{HK}
J.~Hart and K.~Kunen,
Single axioms for odd exponent groups,
\textit{J. Automated Reasoning} \textbf{14} (1995), 383--412.

\bibitem{Kiechle}
H.~Kiechle,
\textit{Theory of K-loops},
Lecture Notes in Mathematics \textbf{1778},
Springer, 2002.

\bibitem{KKP}
M.~K.~Kinyon, K.~Kunen, and J.~D.~Phillips,
A generalization of Moufang and Steiner loops,
\textit{Algebra Universalis} \textbf{48} (2002), 81--101.

\bibitem{Mace4}
W.~W.~McCune,
\textit{Mace4 Reference Manual and Guide},
Tech. Memo ANL/MCS-TM-264, Mathematics and Computer Science Division,
Argonne National Laboratory, Argonne, IL, August 2003.
\texttt{http://www.cs.unm.edu/\symbol{126}mccune/mace4/}

\bibitem{Prover9}
W.~W.~McCune,
\textit{Prover9 Manual},
\texttt{http://www.cs.unm.edu/\symbol{126}mccune/prover9/}

\bibitem{LOOPS}
G.~P.~Nagy and P.~Vojt\v{e}chovsk\'{y}, LOOPS -- a GAP package, version 1.4.0,
Feb. 2007, \texttt{http://www.math.du.edu/loops}

\bibitem{Pflugfelder}
H.~O.~Pflugfelder,
\textit{Quasigroups and Loops: Introduction},
Sigma Series in Pure Math. \textbf{8}, Heldermann Verlag, Berlin, 1990.

\bibitem{PhillipsVojtechovsky}
J.~D.~Phillips and P.~Vojt\v{e}chovsk\'y,
C-loops: An introduction,
\textit{Pub. Math. Debrecen} \textbf{68} (2006), 115--137.

\end{thebibliography}

\end{document}